\documentclass[12pt, a4paper]{article}
\usepackage{mathrsfs, amsmath,amssymb,amsthm,mathtools,setspace,graphicx,array,tikz,url}
\usepackage{cite}
\usepackage[hidelinks]{hyperref}
\usepackage{fancyhdr}
\usepackage{lastpage} % For \pageref{LastPage}
\usepackage[top=3cm, left=3cm, right=3cm, bottom=3cm]{geometry}
\usepackage{enumitem}
\usepackage{titlesec}
\usepackage{etoolbox}

\patchcmd{\thebibliography}{\advance\leftmargin\labelsep}{\labelsep=.5em\leftmargin=2.5em\itemindent=0em\advance\leftmargin\labelsep}

\titlespacing*{\section}{0pt}{3em}{2em}

\newtheorem{conjecture}{Conjecture}[section]
\newtheorem{theorem}{Theorem}[section]
\newtheorem{definition}{Definition}[section]
\newtheorem{example}{Example}[section]

\newtheorem{proposition}{Proposition}[section]
\newtheorem{corollary}{Corollary}[section]
\newtheorem{remark}{Remark}[section]
\newtheorem{namedlemma}{Negation Lemma}[section]

\title{A Note on the Union-closed Sets Conjecture}
\author{
  Mohammad Javad Moghaddas Mehr \\
  \small{Email: m.moghadas11235@gmail.com}
}

\date{September 4, 2023}

\begin{document}
\maketitle
\vspace{5em}
\begin{abstract}
  Let $M$ be a non-zero binary matrix with distinct rows where the rows are closed under certain logical operators. In this article, we investigate the existence of columns containing an equal or greater number of ones than zeros. Specifically, the existence of such columns when the rows of the matrix are closed under \emph{material conditional} leads us to a weaker version of the \emph{Union-Closed Set Conjecture}.
\end{abstract}


\section{Introduction}
\paragraph{} For an integer, $n \in \mathbb{N},$ we define the set $[n]$ as $\{k \in \mathbb{N} \, | \, k \leq n\}.$ A family $\mathscr{F} \subseteq 2^{[n]}$ is called union-closed if $\mathcal{A}, \mathcal{B} \in \mathscr{F}$ implies that $\mathcal{A} \cup \mathcal{B} \in \mathscr{F}.$ In 1979, Péter Frankl \cite{frankl1995} posed a seemingly straightforward conjecture that, over four decades later, still stands as an unsolved puzzle despite substantial efforts and progress dedicated to its understanding. The conjecture states as follows.
\vspace{1em}
\begin{conjecture}\label{conj:re7}
 For every non-empty finite union-closed family of sets, there exists an element that belongs to at least half of the sets in the family.
\end{conjecture}

\paragraph{} Efforts to understand the conjecture have led to notable breakthroughs. We recommend the survey paper by Bruhn and Schaudt \cite{bruhn2015} for details.
Vuckovic and Zivkovic \cite{vuckovic2012} demonstrated the conjecture's validity for any union-closed family $\mathscr{F} \subseteq 2^{[n]}$ where $n \leq 12.$
Roberts and Simpson \cite{roberts2010} established that if $q$ is the minimum cardinality of $\bigcup \mathscr{F}$ taken over all counterexamples $\mathscr{F},$
then any counterexample $\mathscr{F}$ has cardinality at least $4q - 1.$

\paragraph{} In 2017, Ilan Karpas \cite{karpas2017} showed that there exists some absolute constant $c > 0$ such that for any union-closed family $\mathcal{F} \subseteq 2^{[n]},$ if $|\mathcal{F}| \geq \left(\frac{1}{2} - c\right)2^n,$ then there is an element $i \in [n]$ that appears in at least half of the sets in $\mathcal{F}.$

\paragraph{} Recently, Gilmer \cite{gilmer2022} established the first constant lower bound for the conjecture, asserting the existence of an element that belongs to at least 0.01 of the sets in the family. Gilmer claimed that his method could enhance the lower bound to $\left(3 - \sqrt{5}\right)/2 \approx 0.38.$ A few days later, three preprints were published and verified his claim \cite{chase2022, sawin2022, alweiss2022}.

\paragraph{} In this note, we explore binary matrices whose rows are closed under different logical operators. To do so, we need some fundamental definitions.
\vspace{1em}
\begin{definition}
   Let $A = (a_{ij})$ and $B = (b_{ij})$ be two binary matrices of size $n \times m,$
  \begin{enumerate}[parsep=1pt]
    \item Let $\neg A:= (\neg a_{ij}).$
    \item Let $A * B := (a_{ij} * b_{ij})$ for binary operator "$*$".
  \end{enumerate}
\end{definition}

\paragraph{}
Let $M = (m_{ij})$ be a matrix. We denote $M_{i-}$ and $M_{-j}$ as the rows and columns of $M,$ respectively. From this perspective, we can restate Conjecture \ref{conj:re7} as follows.

\begin{conjecture}\label{conj:re6}
   Let $M = (m_{ij})$ be a non-zero binary matrix of size $n \times m$ with distinct rows. If for any arbitrary pair of rows $M_{i-}$ and $M_{j-}$ we have $M_{i-} \lor M_{j-}$ as a row of $M,$ then there exists a column $M_{-k}$ that contains at least $n/2$ ones.
\end{conjecture}

\paragraph {}Finally in Section 7 we will establish a weaker version of Conjecture \ref{conj:re6} by proving the following theorem.
\vspace{1em}
\begin{theorem}
   Let $M = (m_{ij})$ be a non-zero binary matrix of size $n \times m$ with distinct rows. If for any arbitrary pair of rows $M_{i-}$ and $M_{j-}$ we have $\neg M_{i-} \lor M_{j-}$ as a row of $M,$ then there exists a column $M_{-k}$ that contains at least $n/2$ ones.
\end{theorem}

\paragraph {}Let $\mathcal{M}$ be the set of all binary matrices, and denote $\mathcal{N}$ be the set of all finite subsets of natural numbers. Define $\psi : \mathcal{M} \to \mathcal{N}$ as follows

$$
\psi(M_{n \times m}) = \left\{\sum_{i=1}^{n} m_{ij} \, | \, j \in [m]\right\}.
$$
\vspace{1em}
\begin{example}
  
  Let $\mathscr{A}$ be a family of sets as follows 
  $$
  \mathscr{A}=\{\varnothing ,\{1\},\{1,2\},\{2,3,4\},\{1,2,3,4\}\}
  $$

\noindent which is union-closed.

\noindent The matrix representation of $\mathscr{A}$ will be
  $$
  A = \begin{pmatrix}
  0 & 0 & 0 & 0\\
  1 & 0 & 0 & 0\\
  1 & 1 & 0 & 0\\
  0 & 1 & 1 & 1\\
  1 & 1 & 1 & 1
  \end{pmatrix}
  $$ \\
  \noindent then $\psi(A) = \{2,3\},$ so $max(\psi(A)) \geq 5/2.$ Thus, Conjecture \ref{conj:re6} holds in this case.
\end{example}

\paragraph{}
If Conjecture \ref{conj:re6} holds for matrix $A,$ then any row and column permutation of $A$ also satisfies the $\text{Conjecture's}$ conditions.
\vspace{1em}
\begin{definition}
   Matrix $A$ is equivalent to matrix $B,$ denoted as $A \sim B,$ if $B$ can be obtained from $A$ through row and column permutations.
\end{definition}
\paragraph {}
We refer to the pair $(M,*)$ as a space, where $M = (m_{ij})$ is a non-zero binary matrix of size ${n \times m},$ such that its rows are closed under "*" operator. This operator can be either binary or unary.

\section{\texorpdfstring{Negation ($\neg$)}{Negation}}
\paragraph{}
In this section, we will prove our first lemma, which serves as a useful tool and provides valuable insights for dealing with other cases.
\vspace{1em}
\begin{namedlemma}\label{lem:negation}
  For every space $(M_{n \times m},\neg)$ we have $max(\psi(M)) \geq n/2.$
\end{namedlemma}

\begin{proof}
   Without loss of generality, assume that $M_{-j} \neq 0$ for each $j \in [m].$ Define two sets $\mathcal{K}$ and $\mathcal{L}$ as follows
  
    \begin{align*}
      \mathcal{K} &:= \{M_{i-} \,|\, m_{i1}=1\} \\
      \mathcal{L} &:= \{M_{i-} \,|\, m_{i1}=0\}.
    \end{align*}
   
  \noindent If $\mathcal{L} = \varnothing$ there is nothing to prove, so let $\mathcal{L} \neq \varnothing.$ Define matrix $K$ such that its rows come from $\mathcal{K}$ while preserving their order in $M.$ We define matrix $L$ from $\mathcal{L}$ in the same manner, then
  $$
  M \sim \begin{pmatrix}
  K \\
  L \\
  \end{pmatrix}.
  $$

  \noindent For any row $K_{t-}$ there is exactly one row $L_{s-}$ such that $K_{t-} = {\neg}L_{s-}$, and vice versa.

\newpage

  \noindent This establishes a mutual correspondence between the rows of $K$ and the rows of $L$, thus
  $$
  K \sim {\neg}L
  $$

  \noindent which forces $M$ to have an even number of rows and an equal number of ones and zeros in each column.
  So $n=2k$ for some $k \in \mathbb{N}$, then $\psi(M) = \{k\}$. Thus
  $$
  \max(\psi(M)) = k.
  $$
\end{proof}

\section{\texorpdfstring{Alternative denial ($\uparrow$) and Joint denial ($\downarrow$)}{Alternative denial and Joint denial}}

  \paragraph{}
  In this section, we will utilize Lemma \ref{lem:negation} to explore the properties of two binary operators, Alternative Denial and Joint Denial.
  \vspace{1em}
  \begin{proposition}
     For every space $(M_{n \times m},\uparrow)$ we have $\max(\psi(M)) \geq n/2$.
  \end{proposition}
  \begin{proof}
     Since rows of $M$ are closed under "$\uparrow$", then for an arbitrary $M_{i-}$ we have
    $$
    M_{i-} \uparrow M_{i-} = {\neg}M_{i-}
    $$

    \noindent which means rows of $M$ are closed under "$\neg$". By applying the Lemma \ref{lem:negation}, we can conclude
    $\max(\psi(M)) \geq n/2$.
  \end{proof}
  \vspace{1em}
  \begin{proposition}
     For every space $(M_{n \times m},\downarrow)$ we have $\max(\psi(M)) \geq n/2$.
  \end{proposition}
  \begin{proof}
     Since rows of $M$ are closed under "$\downarrow$", then for an arbitrary $M_{i-}$ we have
    $$
    M_{i-} \downarrow M_{i-} = {\neg}M_{i-}
    $$

    \noindent which means rows of $M$ are closed under "$\neg$". By applying the Lemma \ref{lem:negation}, we can conclude
    $\max(\psi(M)) \geq n/2$.
  \end{proof}

\section{\texorpdfstring{Abjunction ($\not\rightarrow$) and Conjunction ($\land$)}{Abjunction and Conjunction}}
  
  \paragraph{}
  In this section, we provide two examples to demonstrate that when rows are closed under \text{Abjunction} or \text{Conjunction}, it does not guarantee the existence of a column with an equal or greater number of ones than zeros.
  \vspace{1em}
  \begin{example}
     Let $A$ be the binary matrix
  
    $$
    A_{(n+1) \times n} = \begin{pmatrix}
    I_{n \times n} \\
    0_{1 \times n}
    \end{pmatrix}
    $$
    
    \noindent then for $n>1$
  
    $$
    \max(\psi(A)) = 1 < (n+1)/2
    $$
  
    \noindent but the rows of matrix $A$ are closed under \text{abjunction} and \text{conjunction}.
  \end{example}

  \paragraph{}
  For an arbitrary $k \in \mathbb{N},$ there exists a matrix $A$ such that its rows are closed under \text{conjunction} and $\max(\psi(A)) = k+1.$
  \vspace{1em}
  \begin{example}
     Let $A$ be the binary matrix
  
    $$
    A_{(n+2) \times (n+1)} =
    \begin{pmatrix}
    \begin{array}{c|c}
    \begin{matrix}
    J_{k \times 1} \\
    0_{(n-k) \times 1} \\
    \end{matrix} & I_n \\
    \hline
    \begin{matrix}
    1 \\
    0 \\
    \end{matrix} & 0_{2 \times n}
    \end{array}
    \end{pmatrix}
    $$
  
    \noindent so $\max(\psi(A)) = k+1.$ For instance, for $k=2$
  
    $$
    A_{7 \times 6} =
    \begin{pmatrix}
  
    \begin{array}{c|c}
  
    \begin{matrix}
    1 \\
    1 \\
    0 \\
    0 \\
    0 \\
    \end{matrix} &
    \begin{matrix}
    1 & 0 & 0 & 0 & 0 \\
    0 & 1 & 0 & 0 & 0 \\
    0 & 0 & 1 & 0 & 0 \\
    0 & 0 & 0 & 1 & 0 \\
    0 & 0 & 0 & 0 & 1 \\
    \end{matrix} \\
    \hline
    \begin{matrix}
    1 \\
    0 \\
    \end{matrix} &
    \begin{matrix}
    0 & 0 & 0 & 0 & 0 \\
    0 & 0 & 0 & 0 & 0 \\
    \end{matrix}
  
    \end{array}

    \end{pmatrix}
    $$
  
    \noindent thus
  
    $$
    \max(\psi(M)) = 3 < 7/2.
    $$
  \end{example}

\section{\texorpdfstring{Biconditional ($\leftrightarrow$) and Exclusive or ($\not\leftrightarrow$)}{Biconditional and Exclusive or}}
  
\paragraph{}
In this section, we will explore two logical operators that form an additive group on rows of a binary matrix.
\vspace{1em}
\begin{proposition}
   For every space $(M_{n \times m},\leftrightarrow)$ we have $\max(\psi(M)) \geq n/2.$
\end{proposition}

\begin{proof}
   Define matrices $K$ and $L,$ as in the proof of the Lemma \ref{lem:negation}, so
\newpage

  $$
  M \sim \begin{pmatrix}
  K \\
  L \\
  \end{pmatrix}.
  $$

\noindent Define the function $g$ as follows

\begin{align*}
  g: \mathcal{L} &\to \mathcal{K} \\
  g(L_{t-}) &= L_{1-}+L_{t-}.
\end{align*}

  \noindent If $g(A) = g(B)$ then $L_{1-} + A = L_{1-} + B.$ Since $\{M_{i-}\}_{i \in [n]}$ is an additive group, then $A = B.$ So $g$ is injective, and as a result,

  $$
  |\mathcal{L}| \leq |\mathcal{K}|.
  $$

  \noindent Thus

  \begin{align*}
  \max(\psi(M)) &\geq (|\mathcal{K}|+|\mathcal{L}|)/2 \\
  &= n/2.
  \end{align*}

\end{proof}
\vspace{1em}
\begin{proposition}
   For every space $(M_{n \times m},\not\leftrightarrow)$ we have $\max(\psi(M)) \geq n/2.$
\end{proposition}

\begin{proof}
   Define matrices $K$ and $L,$ as in the proof of the Lemma \ref{lem:negation} , so

  $$
  M \sim \begin{pmatrix}
  K \\
  L \\
  \end{pmatrix}.
  $$

  \noindent Define the function $g$ as follows

  \begin{align*}
    g: \mathcal{L} &\to \mathcal{K} \\
    g(L_{t-}) &= K_{1-}+L_{t-}.
  \end{align*}

  \noindent If $g(A) = g(B)$ then $K_{1-} + A = K_{1-} + B.$ Since $\{M_{i-}\}_{i \in [n]}$ is an additive group, then $A = B.$ So $g$ is injective, and as a result,

  $$
  |\mathcal{L}| \leq |\mathcal{K}|.
  $$

  \noindent Thus

    \begin{align*}
      \max(\psi(M)) &\geq (|\mathcal{K}|+|\mathcal{L}|)/2 \\
       &= n/2.
    \end{align*}

\end{proof}

\section{\texorpdfstring{Topology ($\cup, \cap$)}{Topology}}
\paragraph{}
If $\mathscr{F} \subseteq 2^{[n]}$ is a family of sets that is closed under "$\cup$" and "$\cap$", it
implies that we may define a topology on $\bigcup \mathscr{F}.$ However, being closed under "$\cup$"
and "$\cap$" doesn't guarantee that $\varnothing \in \mathscr{F}.$ Nevertheless, if $\bigcap \mathscr{F} \neq \varnothing,$
then there exists an $a \in \bigcup \mathscr{F}$ such that it belongs to each member of $\mathscr{F}.$

\begin{theorem}
     Let $(X, \tau)$ be a finite topological space, then there exists an $x_0 \in X$ that appears in at least half of the members of $\tau$.
\end{theorem}

\begin{proof}
     Without loss of generality, assume $X = [n]$ for some $n \in \mathbb{N}$. Let $\mathscr{X}$ be a family of members of $\tau$ with the minimum cardinality among the members of $\tau$, excluding the empty set. We select the smallest set $\mathcal{B} \in \mathscr{X}$ according to the lexicographic order. It is guaranteed that for $\mathcal{A} \in \tau$, we have either $\mathcal{B} \subseteq \mathcal{A}$ or $\mathcal{B} \cap \mathcal{A} = \varnothing$. Define sets $\mathcal{K}$ and $\mathcal{L}$ as follows

    \begin{align*}
        \mathcal{K} &= \{\mathcal{A} \in \tau \,|\, \mathcal{B} \subseteq \mathcal{A}\} \\
        \mathcal{L} &= \{\mathcal{A} \in \tau \,|\, \mathcal{B} \cap \mathcal{A} = \varnothing\}.
    \end{align*}

    \noindent If $\mathcal{L} = \varnothing$, there is nothing to prove, so let $\mathcal{L} \neq \varnothing$. Define the function $g$ as follows

    \begin{align*}
      g: \mathcal{L} &\to \mathcal{K} \\
      \quad g(\mathcal{A}) &= \mathcal{B} \cup \mathcal{A}.
    \end{align*}

    \noindent If $g(\mathcal{A}) = g(\mathcal{C})$, then $\mathcal{A} \cup \mathcal{B} = \mathcal{C} \cup \mathcal{B}$. Since $\mathcal{B}$ has an empty intersection with $\mathcal{A}$ and $\mathcal{C}$, then $\mathcal{A} = \mathcal{C}$. So $g$ is injective, and as a result, $|\mathcal{L}| \leq |\mathcal{K}|$, thus members of $\mathcal{B}$ appear in at least half of the members of $\tau$.

\end{proof}
\vspace{1em}
\begin{corollary}
     Let $M$ be a non-zero binary matrix of size $n \times m$, where its rows are closed under "$\land$" and "$\lor$", then $\max(\psi(M)) \geq \frac{n}{2}$.
\end{corollary}

\section{\texorpdfstring{Material Conditional ($\rightarrow$)}{Material Conditional}}

\paragraph{}
In this section, we will establish our main result. Before proceeding, we need some tools that will help us
move forward more easily.

\paragraph{}
Let $M_{n \times m}$ be a non-zero binary matrix. We define $\widetilde{M}_{n \times m}={\neg}M.$
If rows of $M$ are closed under the binary operator "*", then rows of $\widetilde{M}$ would also be closed
under the binary
operator "$\widetilde{*}$" defined as follows

$$
\widetilde{*}(\widetilde{A},\widetilde{B}) := {\neg}*(A,B)
$$

\noindent where $\widetilde{A}, \widetilde{B} \in \{\widetilde{M}_{i-}\}_{i \in [n]}$ and $A, B \in \{M_{i-}\}_{i \in
[n]}.$

\begin{example}\label{ex:re1}

  Let $(M_{n \times m},*)$ be a non-zero space such that
  for $A, B \in \{M_{i-}\}_{i \in [n]},$ $*(A,B) = {\neg}A \lor B.$ Then
  for $\widetilde{A}, \widetilde{B} \in \{\widetilde{M}_{i-}\}_{i \in [n]}$
  $$
  \widetilde{*}(\widetilde{A},\widetilde{B}) = {\neg}*(A,B) = {\neg}({\neg}A \lor B) = {\neg}\widetilde{A}
  \land
  \widetilde{B}.
  $$

\end{example}
\vspace{1em}
\begin{definition}
     Let $\mathscr{F}$ be a non-empty finite family of sets. A \emph{basis} for the family $\mathscr{F}$ is a disjoint collection $\mathscr{B} = \{\mathcal{B}_i\}_{i \in \gamma}$ of members of $\mathscr{F}$ such that for every member $\mathcal{A}$ in $\mathscr{F}$, there exists an index set $\alpha \subseteq \gamma$ such that $\mathcal{A} = \bigcup_{i \in \alpha} \mathcal{B}_i$.
\end{definition}

\begin{definition}
     Let $M_{n \times m}$ be a non-zero binary matrix. A \emph{basis} for the matrix $M$ is a collection $\mathcal{V} = \{v_i\}_{i \in \gamma}$ of rows of $M$, where $v_t$ and $v_k$ are orthogonal for $t \neq k$, such that for every row $A$ of $M$, there exists an index set $\alpha \subseteq \gamma$ such that $A = \bigvee_{i \in \alpha} v_i$.
\end{definition}

\begin{proposition}
  Let $\mathscr{F}$ be a finite family of sets. If $\mathscr{F}$ has a basis, then it is \emph{unique}.
\end{proposition}
\begin{proof}
     Suppose $\mathscr{F}$ has two different bases $\mathscr{B} = \{\mathcal{B}_i\}_{i \in \gamma}$ and $\mathscr{C} = \{\mathcal{C}_i\}_{i \in \lambda}$. Without loss of generality, let $\mathcal{C}_t \in \mathscr{C}$ such that $\mathcal{C}_t \notin \mathscr{B}$, so
    $$
    \mathcal{C}_t = \bigcup_{j \in \alpha} \mathcal{B}_j,
    $$
    where $\mathcal{B}_j \neq \mathcal{C}_t$ and

  $$
  \mathcal{B}_j = \bigcup_{i \in \beta_j} \mathcal{C}_{i}
  $$

  \noindent where $\mathcal{C}_i \neq \mathcal{C}_t.$ Then

  $$
  \mathcal{C}_t = \bigcup_{j \in \alpha} \bigcup_{i \in \beta_j} \mathcal{C}_{i} = \bigcup_{i \in \boldsymbol{\cup}_{j \in \alpha} \beta_j}\mathcal{C}_{i}
  $$

  \noindent where $\mathcal{C}_i \neq \mathcal{C}_t.$ So $\mathcal{C}_t$ is represented by other members of $\mathscr{C}$, which is a contradiction.
  
\end{proof}

\vspace{1em}
\begin{proposition}
   Let $\mathscr{F}$ be a non-empty finite family of sets such that for $\mathcal{A},\mathcal{D} \in \mathscr{F},$ $\mathcal{A} - \mathcal{D} \in \mathscr{F}$ and $\mathcal{A} \cap \mathcal{D} \in \mathscr{F}.$ Then $\mathscr{F}$ has a \emph{unique} basis $\mathscr{B}=\{\mathcal{B}_i\}_{i \in \gamma},$ where for every $\mathcal{A} \in \mathscr{F}$ there is a unique index set $\alpha \subseteq \gamma$ such that $\mathcal{A}=\bigcup_{i \in \alpha} \mathcal{B}_i.$ So $\mathcal{A}$ is uniquely determined by members of $\mathscr{B}.$
\end{proposition}

\begin{proof}
   We construct $\mathscr{B}$ inductively. Start with $\mathcal{F}_0 = \mathscr{F}.$ In each step, let $\mathcal{F}_{i+1}$ be obtained from $\mathcal{F}_{i}$ by removing an element $\mathcal{C}$ that can be expressed as the union of other elements of $\mathcal{F}_{i}.$ This process terminates in less than $|\mathscr{F}|$ steps. Let $\mathscr{B}=\{\mathcal{B}_i\}_{i \in \gamma}$, which is obtained from the last step, so the members of $\mathscr{F}$ can be represented as the union of members of $\mathscr{B}.$

   \bigbreak
   \noindent Assume there is a distinct pair $\mathcal{B}_t, \mathcal{B}_k \in \mathscr{B}$ such that $\mathcal{B}_t \cap \mathcal{B}_k \neq \varnothing.$ Without loss of generality, we can assume $\mathcal{B}_t-\mathcal{B}_k \neq \varnothing.$ Define $\mathcal{A}_1=\mathcal{B}_t-\mathcal{B}_k$ and $\mathcal{A}_2=\mathcal{B}_t \cap \mathcal{B}_k,$ obviously $\mathcal{A}_1$ and $\mathcal{A}_2$ are proper subsets of $\mathcal{B}_t$ and $\mathcal{A}_1, \mathcal{A}_2 \in \mathscr{F},$ hence they can be represented by members of $\mathscr{B},$ so $\mathcal{A}_1 = \bigcup_{i \in \alpha} \mathcal{B}_{i},$ $ \mathcal{A}_2 = \bigcup_{i \in \beta} \mathcal{B}_{i}$ where $\mathcal{B}_i \neq \mathcal{B}_t$ for $i \in \alpha \cup \beta.$ Thus
  $$
  \mathcal{B}_t = \mathcal{A}_1 \cup \mathcal{A}_2
  = \bigcup_{i \in \alpha \cup \beta} \mathcal{B}_{i}
  $$

  \noindent where $\mathcal{B}_i \neq \mathcal{B}_t.$ So $\mathcal{B}_t$ is represented by members of $\mathscr{B}$ which is a contradiction, so different members of $\mathscr{B}$ have empty intersection.
  \bigbreak
  
  \noindent It remains to show that every $\mathcal{A} \in \mathscr{F}$ is uniquely determined by members of $\mathscr{B} = \{\mathcal{B}_i\}_{i \in \gamma}.$ Assume $\mathcal{A}=\bigcup_{i \in \alpha}\mathcal{B}_i$ and $\mathcal{A} = \bigcup_{j \in \beta} \mathcal{B}_j$ where $\alpha, \beta \subseteq \gamma$ and $\alpha \neq \beta.$
  Let $\mathcal{B}_k \in \{\mathcal{B}_i\}_{i \in \alpha}$ thus $\mathcal{B}_k \subseteq \bigcup_{j \in \beta} \mathcal{B}_j$ so there is $\mathcal{B}_t \in \{\mathcal{B}_j\}_{j \in \beta}$ such that $\mathcal{B}_k \cap \mathcal{B}_t \neq \varnothing$ which is a contradiction. So $\mathcal{A}$ is uniquely determined by members of $\mathscr{B}.$
  
\end{proof}

\vspace{1em}
\begin{corollary}\label{cor:re2}

 Let $M_{n \times m}$ be a non-zero binary matrix such that for $A, B \in \{M_{i-}\}_{i \in [n]},$ $A \land \neg B \in \{M_{i-}\}_{i \in [n]}$ and $A \land B \in \{M_{i-}\}_{i \in [n]}.$ Then $M$ has a \textit{unique} basis $\mathcal{V}=\{v_i\}_{i \in \gamma},$ where for every $A \in \{M_{i-}\}_{i \in [n]}$ there is a unique index set $\alpha \subseteq \gamma$ such that $A= \bigvee_{i \in \alpha} v_i.$ So $A$ is uniquely determined by members of $\mathcal{V}.$

\end{corollary}

\paragraph{}
In the following proposition, we will demonstrate the necessary conditions for the matrix $\widetilde{M},$ corresponding to the space $(M_{n \times m},\rightarrow),$ to have a basis.

\bigbreak
\begin{proposition}\label{prop:re3}

 Let $(M_{n \times m},\rightarrow)$ be a non-zero space. Then
\begin{enumerate}[parsep=1pt]
    \item For $\widetilde{A}, \widetilde{B} \in \{\widetilde{M}_{i-}\}_{i \in [n]}$ we have $\widetilde{A} \land \neg \widetilde{B} \in \{\widetilde{M}_{i-}\}_{i \in [n]}.$
    \item For $\widetilde{A}, \widetilde{B} \in \{\widetilde{M}_{i-}\}_{i \in [n]}$ we have $\widetilde{A} \land \widetilde{B} \in \{\widetilde{M}_{i-}\}_{i \in [n]}.$
\end{enumerate}

\end{proposition}

\begin{proof}

 We conclude $(1)$ from Example \ref{ex:re1}.

$(2)$ Let $\widetilde{A}, \widetilde{B} \in \{\widetilde{M}_{i-}\}_{i \in [n]}.$ By $(1)$ $\widetilde{A} \land \neg (\widetilde{A} \land \neg \widetilde{B}) \in \{\widetilde{M}_{i-}\}_{i \in [n]}.$ Since

\begin{align*}
    \widetilde{A} \land \neg (\widetilde{A} \land \neg \widetilde{B}) &= \widetilde{A} \land ({\neg}\widetilde{A} \lor \widetilde{B}) \\
    &= \widetilde{A} \land \widetilde{B}
\end{align*}

\noindent then $\widetilde{A} \land \widetilde{B} \in \{\widetilde{M}_{i-}\}_{i \in [n]}.$

\end{proof}

\vspace{1em}
\begin{remark}

 Let $(M_{n \times m},\rightarrow)$ be a non-zero space. Then $\widetilde{M}$ has a basis $\mathcal{V}=\{v_i\}_{i \in \gamma}$ by Proposition \ref{prop:re3} and Corollary \ref{cor:re2}.

\end{remark}

\newpage

\begin{remark}\label{rem:re4}

 Suppose $M = (m_{ij})$ be a non-zero binary matrix of size $n \times m$ and $k \in [m].$ Clearly $\sum_{i=1}^{n} m_{ik} \geq n/2$ if and only if $\sum_{i=1}^{n} \widetilde{m}_{ik} \leq n/2.$

\end{remark}

\begin{theorem}\label{thm:re5}

 For every space $(M,\rightarrow)$ we have $\max(\psi(M)) \geq n/2.$

\end{theorem}

\begin{proof}

Let $\mathcal{V}=\{v_i\}_{i \in \gamma}$ be a basis for $\widetilde{M}.$ We define sets $\mathcal{K}$ and $\mathcal{L}$ as follows

  \begin{align*}
    \mathcal{K} &:= \{\widetilde{M}_{i-}|\text{ for } \widetilde{M}_{i-}=\bigvee_{i \in \alpha} v_i \text{, where } 1 \in \alpha\} \\
    \mathcal{L} &:= \{\widetilde{M}_{j-}|\text{ for } \widetilde{M}_{j-}=\bigvee_{j \in \beta} v_j \text{, where } 1 \notin \beta\}.
  \end{align*}

\noindent
Clearly $\mathcal{K}$ and $\mathcal{L}$ are well defined, due to the unique representation of rows of $\widetilde{M}$ by members of $\mathcal{V}.$ Define matrix $K$ where its rows come from $\mathcal{K}$ while preserving their order in $\widetilde{M}.$

\noindent We define matrix $L$ from $\mathcal{L}$ in the same manner, then

$$
\widetilde{M} \sim \begin{pmatrix}
K \\
L \\
\end{pmatrix}.
$$

\noindent Define function $g$ as follows

    \begin{align*}
      g: \mathcal{K} &\to \mathcal{L} \\
      g(\widetilde{M}_{k-}) &= \widetilde{M}_{k-} \land \neg v_1.
    \end{align*}

\noindent Since $v_t \land v_k=0$ for any two distinct members $v_t, v_k \in \mathcal{V}$ and every row of $\widetilde{M}$ has a unique representation, we conclude that $g$ is an injective function. So $|\mathcal{K}| \leq |\mathcal{L}|$ which means $v_1$ appears in equal or less than half of rows of $\widetilde{M}.$ Thus, there is a column $\widetilde{M}_{-t} \in \{\widetilde{M}_{-j}\}_{j \in [m]}$ such that $\sum_{i=1}^{n} \widetilde{m}_{it} \leq n/2.$ So by Remark \ref{rem:re4}

$$
\sum_{i=1}^{n} m_{it} \geq n/2
$$

\noindent which implies that

$$
\max(\psi(M))\geq n/2.
$$

\end{proof}
   
\paragraph{}
We end this paper by mentioning that Theorem \ref{thm:re5} provides a weaker version of Conjecture \ref{conj:re6}. To see this, consider a binary matrix $M$ where ${\neg}A \lor B \in \{M_{i-}\}_{i \in [n]}$ for $A, B \in \{M_{i-}\}_{i \in [n]}.$ By Proposition \ref{prop:re3}, we also know that ${\neg}A \land {\neg}B \in \{\widetilde{M}_{i-}\}_{i \in [n]}.$ So we can conclude that $A \lor B \in \{M_{i-}\}_{i \in [n]}.$

\end{document}